\documentclass[10pt]{article}

\usepackage[dvips]{graphicx}
\usepackage{amsmath,amsfonts,amssymb}
\usepackage{a4wide}

\newtheorem{Proposition}{Proposition}[section]
\newtheorem{Theorem}{Theorem}[section]
\newtheorem{Definition}{Definition}[section]
\newtheorem{Remark}{Remark}[section]
\newtheorem{Corollary}{Corollary}[section]
\newtheorem{Lemma}{Lemma}[section]
\newtheorem{Example}{Example}[section]
\newtheorem{Examples}{Examples}[section]

\def\N{{\mathbb N}}
\def\Z{{\mathbb Z}}

\def\C{{\mathbb C}}

\def \deg {{\rm deg\,}}

\def \lra {\longrightarrow}

\def \P {\mathbb P}

\def \N {{\bf N}}
\def \M {${\cal M}_\Gamma$}

\def \I {{\mathfrak I}}
\def \M {{\mathfrak M}}
\def \J {{\mathfrak J}}

\def \basisI {{\cal B}(\I)}
\def \basis {{\cal B}}
\def \ffdd {\vrule height 4pt width 3pt depth 2pt}
\def \finedim {\par \nobreak \rightline{\ffdd \qquad}}
\def \l {\lambda}
\def \a {\alpha}
\def \b {\beta}

\begin{document}

\title{ \textbf{Strongly inessential elements of a perfect height 2 ideal}}

\author{\\ Giannina Beccari, Carla Massaza \\
Dipartimento di Scienze Matematiche,  Politecnico Torino \\
Corso Duca degli Abruzzi, 24 - 10129 Torino - Italy\\
 \\}

\date{12/12/2011}
\maketitle

\begin{abstract}
In this paper we expand on some results exposed in  a previous one, in which we introduced the concept of inessential and strongly inessential generators in a standard basis of a saturated homogeneous ideal. The appearance  of strongly inessential elements seemed to be a non generic situation; in this paper we analyze their presence in a perfect height 2 ideal with the greatest number of generators, according to Dubreil's inequality.

\end{abstract}

\section{Introduction}

In a previous paper \cite{BM2} we introduced the concept of {\it strongly inessential element} (briefly s.i.) in a homogeneous ideal $\I \subset K[X_1,...,X_n]$. Our first idea, when we started to think about essential and inessential elements of a standard basis (see \cite{BM2}, n.3), was that every homogeneous ideal should have a standard basis consisting of essential forms, but we very soon found  many counterexamples. Therefore, our next conjecture was that the assertion might be true for a sufficiently general ideal. In this paper  we thus  investigate the structure of e-maximal bases  (\cite{BM2} Definition 5.1) and, as a consequence ( \cite{BM2} Theorem 5.1), the presence of s.i. elements, in what seemed to be the easiest situation, that is  when $\I$ is a perfect height 2 ideal. In this case, it is possible to associate to every $\basisI$  a Hilbert-Burch matrix( \cite{N},\cite{St}) and to decide the nature of the forms of $\basisI$, with respect to essentiality, just looking at the ideals generated by the entries of its columns (\cite{BM2}).

 We observe that, if the multiplicity $e(\I)$ (\cite{F},\cite{Ha}, \cite{ZS})\ is low, our first idea was correct; more precisely, if $e(\I)<6$, then every standard basis consists of essential elements, while, if $6< e(\I) <9$, $\I$ has at least a standard basis whose elements are all essential.

  To deal with the problem when the multiplicity is $\geq 9$, we observe that  strong inessentiality  is preserved modulo a regular sequence (while essentiality is not). So, the first case to be considered seems to be the one of zero depth. As the general case still appears hard to be analyzed, we replace the family of all perfect height 2 ideals with its subfamily  $F = \bigcup_{n\geq 2} F[n]$, where $F[n]$ is the set of all perfect height two ideals in  $S = K[X_1,...,X_n], n\geq 2$, whose standard bases  are of maximal cardinality with respect to Dubreil's inequality (\cite{D}). In a previous paper  \cite{BM1}, in fact, we found a description of $F$ that is of help in dealing with the problem considered here. So, as we restrict our attention to the ideals of zero depth, we study $F[2]$. For every ideal $\I \in F[2]$, we produce a {\it canonical} Hilbert matrix, with the property that its corresponding basis is e-maximal, which means that its inessential elements are  s.i.. Using such a matrix, we prove that the number of s.i. elements appearing in an e-maximal basis is linked to the greatest common divisor $\Phi$ of its generators of minimal degree $\alpha(\I)$; in fact, it depends on the decomposition of $\Phi$ into linear factors  (see Theorem \ref{Th 3.1}). More precisely, we prove that $\I$ has a basis of essential elements iff all the linear factors of $\Phi$ are distinct; therefore, the generic $\I\in F[2]$ has this a property.

  The description of the e-maximal bases is much more complicated when we pass from $F[2]$ to $F[3]$.
 The Hilbert-Burch matrix of any element $\I\in F[3]$ can be obtained by lifting the Hilbert-Burch matrix of its image $\bar{\I} \in F[2]$ modulo any linear form, regular for $S/\I$ (\cite{BM1}); however, it may happen that there exists some $\widetilde{\I}$ with the same number of s.i. elements of $\bar{\I}$ in any e-maximal basis, among the ideals of $F[3]$ lifting $\bar{\I} \in F[2]$, but there are also cases in which  no lifting of $\bar{\I}$ preserves a s.i. element. We prove that the greatest expected number of s.i. generators in a standard basis of $\I \in F[3]$ is $\a(\I)-2$ and that this number is attained. So, we focus on the set ${\cal S}\subset F[3]$ of the ideals with $\a(\I)-2$  s.i. generators in their e-maximal bases, finding some of their properties and giving examples. In particular, we  completely describe the ideals $\I$ generated in two different degrees, with $\a(\I) =3$ and a s.i. element in any e-maximal basis.

\section{Background and Notation}

Let $S=K[X_1,...,X_n]$, where $K$ is an algebraically closed field, be the coordinate ring of $\P^{n-1}$, $\I =\bigoplus \I_d, d\in \N$, a homogeneous ideal of $S$, and $\M = (X_1,...,X_n)$ be the irrelevant maximal ideal. We recall some basic definitions.

The Hilbert function of $S/\I$ \ (\cite{Ma}), denoted $ H(S/\I,-)$, is the function defined by $$ H(S/\I,t) = dim_K (S/\I)_t.$$ It is well known that for $t\gg 0$ the function $H(S/\I,t)$ is a polynomial, with rational coefficients, of degree $ r(S/\I)-1$, where $r(S/\I)$ is the Krull dimension of $S/\I$.

If $\Delta$ denotes the difference operator on maps from $\Z$ to $\Z$, defined by  $\Delta \phi (t) = \phi(t) - \phi(t-1)$, the function $$\Gamma(\I,t) = \Delta^{r(S/\I)}H(S/\I,t)$$  is called the Castelnuovo function of $\I$,  while $\Delta^{r(S/\I)}H(S/\I,t)$ is, for large $t$, a natural number $e(\I)$, independent on $t$, which is called the {\it multiplicity} of $S/\I$, or also of $\I$.

\begin{Definition} \label{Def. 1.1}{\rm  (\cite{DGM}) \  A  {\it standard basis} $\basisI$ \ of $\I$ is an ordered set of forms of $S$, generating $\I$, such that its elements of degree $d$ define a $K$-basis of $\I_d/(\I_{d-1}S_1)$(\cite {Da},\cite{DGM},\cite{Ca}, \cite{KR}).}
\end{Definition}

It is well known (\cite{DGM}) that the degree vector of $\basisI$ , with non decreasing entries, does not depend on the basis; $\alpha(\I)$ denotes its first entry, $\nu(\I)$ the number of entries, $\nu(\I,t)$ the number of entries equal to $t$. Moreover, if $ht(\I) > 1, \beta(\I)$ is the minimal degree $t$ such that $GCD(\I_t)=1$.

The following theorem links $\alpha(\I)$ to $\nu(\I)$.

\begin{Theorem}\label{ Th 1.1} {\rm (Dubreil)(\cite{Da, D, DGM}) \  Let $\I$ be a homogeneous perfect height 2 ideal. Then $\nu(\I) \leq \alpha(\I) + 1$. }

\end{Theorem}

According to \cite{BM1}, ${\cal F}[n]$ denotes the set of all the homogeneous perfect height 2 ideals of $ S = K[X_1,...,X_n]$  such that $\nu(\I) = \alpha(\I) + 1$; in this paper they are called �Dubreil' s ideals�. In the special case $n=2$, Theorem 1.7 ii) of {\cite {BM1}} gives a description of every ideal of ${\cal F}[2]$ involving the greatest common divisor $\Phi$ of its elements of degree $\alpha(\I)$ and a decomposition of $\Phi$ as a product of forms.

\bigskip

A refinement of Theorem \ref{ Th 1.1} (\cite{Ca}) says, in particular, that, for every perfect height 2 ideal $\I$ in $S$
\begin{equation}
  t\leq \beta(\I)\Rightarrow \nu(\I,t) \leq -\Delta \Gamma (\I,t). \label{eq 1.1}
\end{equation}

We say that $\nu(\I,t)$ is maximal when equality holds in (\ref{eq 1.1}).

\bigskip
 If $\I$ is a perfect height 2 ideal, then a minimal resolution of $S/\I$ is defined by a Hilbert-Burch (shortly H.B.) matrix $M(\I)$ which, in turn, is uniquely determined by a standard basis $\basisI$ \ and by a minimal basis of its module of syzygies \ $Syz \, {\cal B} (\I)$. Its corresponding degree matrix $\partial M(\I)$ is uniquely determined by $\I$.

 \bigskip
 We need some results, widely explained in \cite{BDM1, BDM2}, that we summarize as follows.

 \begin{Theorem} \label{Th 1.2} {\rm Let $\I$ be a perfect height $2$ ideal, $ p+1 $ a degree in which the number $\nu(\I,p+1)$ of generators  in degree $ p+1 $ satisfies the following relation of maximality with respect to Dubreil-Campanella inequality

\begin{equation}
\nu(\I,p+1)= \Gamma(\I,p)- \Gamma(\I,p+1), \label{eq 2.1}
\end{equation}
 $D$ the greatest common divisor of $\I_p$. Then $\I$ admits a basis  $\basis = (DF_1,...,DF_m,G_1,...,G_n)$, where $(DF_1,...,DF_m)S = \I_p $, so that $\I$ splits into two ideals  $\I' =(F_1,...,F_m)$  and  $\I'' = (D, G_1,...,G_n)$, which are still perfect of height 2. Moreover, there is a H.B. matrix $M(\I)$ with respect to $\basis$, with the following shape

$$M(\I) = \begin{pmatrix}A  &  0  \\
                          B  &  C  \end{pmatrix},$$
where

\begin{description}
\item i) {$A \in K^{(m-1)\times m}$ is a H.B. matrix of $\I'$,}
\item ii) {A H.B. matrix of $\I''$ is $(B''\ C)$, where $B'' = B \ ^t(F_1...F_m)$,}
\item iii) {det $C = D$}
\end{description}}
\end{Theorem}

 \section{Strongly inessential elements of an ideal: recalls and complements}

Let $\I= \oplus\I_d, d\in \N, \I_d \subset S_d$  be a homogeneous ideal of $ S = K[X_1,...,X_n]$. We recall some definitions and results appearing in  \cite{BM2}.

\begin{Definition} \label{Def 2.1} {\rm (\cite{BM2}  Definition 3.1)  An element $f$ of a standard basis $\basisI$  is called an {\it inessential} generator of $\I$ with respect to $\basisI$ iff  $$ \exists t \in    \N, fM^t \subseteq (\basisI-\{f\})S. $$  Otherwise we say that $f$ is an {\it essential} generator of $\I$ with respect to $\basisI$.}
\end{Definition}
\bigskip

In the special case of perfect height 2 ideals, the essentiality of the r-th element $f_r$ of $\basisI$ can be read on the ideal $\I_{C_r}$ generated by the entries of the r-th column of any matrix of $Syz \,\basisI$.  In fact Proposition 4.1 of \cite{BM2} says what follows.

\begin{Proposition} \label{Prop 2.2} {\rm Let $\I$ be a perfect codimension 2 ideal of $S$. Then $f_r \in \basisI$ is inessential for $\basisI$ iff the condition

$$ (\exists t \in \N )  \ M^t \subseteq \I_{C_r}$$
is satisfied.}

\end{Proposition}

\begin{Definition} \label{Def 2.3} {\rm (\cite{BM2} Definition 3.2) An element $f \in \I_d$ is strongly inessential (s.i.) iff $ f\notin (\I_{d-1}) S$ and it is inessential with respect to any standard basis containing it.}
\end{Definition}

\begin{Definition} \label{ Def 2.4} {\rm (\cite{BM2} Definition 5.1) A standard basis is called e-maximal iff it has, in every degree $d$, exactly $\nu_e(d)$ essential generators, where $\nu_e(d)$ is the greatest number of essential generatores of degree $d$ appearing in a standard basis of $\I$.}

\end{Definition}

\begin{Theorem} \label{Th 2.5} {\rm (\cite{BM2} Theorem 5.1) A standard basis is e-maximal iff its inessential elements are strongly inessential.}

\end{Theorem}

\bigskip

Starting from Theorem \ref{Th 2.5} we can prove the following statement.

\begin{Proposition} \label{Prop 2.6} {\rm The ideal $\I\subset S$ admits a basis of essential elements iff none of its elements is s.i..}

\end{Proposition}

\medskip \noindent
{\bf Proof.} Proposition 5.2 of \cite{BM2} says that two different e-maximal bases contain the same number of inessential elements. So, $\I$ has a basis of essential elements iff all its e-maximal bases do not contain inessential elements, and we know that they should be s.i., thanks to Theorem \ref{Th 2.5}. Now, every s.i. element can be considered as an entry of a standard basis $\basisI$ and from any standard basis $\basisI$ it is possible to produce an e-maximal basis ${\cal B}_M(\I)$, containing as a subset all the s.i. elements appearing in $\basisI$ (see Proposition 5.4 in \cite{BM2}). So, the e-maximal bases do not contain inessential elements iff s.i. elements do not exist in $\I$.

\finedim

 In other words, $\I$ admits a basis of essential elements iff one of its  e-maximal basis (and, as a consequence, all of them) consists of essential elements and this is equivalent to say that $\I$ does not contain s.i. forms.

\bigskip
Next proposition says that a s.i. element of $\I$ preserves its property modulo a linear form, regular for $S/\I$.

We will use the following notation.

{\bf Notation}  If $z$ is any element of $S=K[X_1,...,X_n]$ and $\phi: S\lra S/zS$ is the canonical morphism, then we set : $ \phi(s)= \bar{s}, \forall  s\in S$  and $\phi({\cal A}) =\bar{\cal{A}}$ for any subset ${\cal A} \subseteq S$, if the element $z$ can be understood.

We need the following lemma.

\begin{Lemma}\label{Lemma 2.7} {\rm \cite{Da, DGM} If $\basis$ is a standard basis of $\I$ and $z\in S$ is a linear form, regular for $S/\I$, then $\bar{ \basis}$ is a standard basis of $\bar{\I}$.}
\end{Lemma}

\medskip

\begin{Proposition} \label{Prop 2.7} {\rm Let $s\in \I$ be a s.i. element and $z$ a linear form regular for $S/\I$. Then $\bar{s} \in \bar{\I}$ is s.i..}
\end{Proposition}

\medskip \noindent
{\bf Proof.} Without any loss of generality we can suppose $z=X_1$. At first we notice that if $s$ is inessential for $\basisI = \basis$, then $\bar{s}$ is inessential for the standard basis $\bar{\cal B} (\bar{\I})= \bar{\basis}$ of $\bar{\I}$. In fact we have:
$$s \ \M^t \subseteq (\basis -\{s\})S \Rightarrow \bar{s} \ \bar{\M}^t \subseteq ({\cal\bar{B}}-\{ \bar{s}\}).$$
Now, let us suppose $s$ to be s.i. and consider a standard basis $\basis$ containing it, say  $\basis =(b_1,b_2,...,b_h)$, where $b_i=s$. Then $\bar{\basis}$ is a standard basis of $\bar{\I}$, containing $\bar{b}_i= \bar{s}$ and any other standard basis $\cal{C}$ of $\bar{\I}$ is of the form  ${\cal C} =\bar{\basis} P$, where $P=(p_{ji})$ is an invertible matrix, whose entries are forms in $K[X_2,...,X_n]$. Let us observe that $\bar{s}= \bar{b}_i \in  {\cal C}$ iff $p_{ii}\not= 0$ and $p_{ij}=0$ when $j\not=i$. As a consequence, ${\cal B}' = \basis P$ is a standard basis containing \ $s= b_i$; in $\basis'$ the element $s$ is inessential, as it is so in every basis in which it appears. The first part of the proof allows to conclude that $\bar{s}$ is inessential for ${\cal C}$.

\finedim

In section 5 we will see that the lifting of a s.i. element of $\bar{\I}$ is not necessarily s.i. in $\I$. (see Remark \ref{Rem 4.1}).

A consequence of Proposition \ref{Prop 2.7} is that if the image $\bar{\I}$ of $\I\subset K[X_1,...,X_n]$ modulo a maximal regular sequence does not contain any s.i. element, the same property holds for $\I$.  So, it seems convenient to start considering the problem of the presence of s.i. elements when $depth \ (S/\I)=0$ (see section 3).

\bigskip
In the sequel we use the following statement (see Theorem \ref{Th 1.2} for notation).

\begin{Theorem} \label{Th 2.8} { \rm Let $\I\subset S$ be a perfect height $ 2 $ ideal and $ p+1 $ a degree in which the maximality condition (\ref{eq 2.1}) is verified. The following statements hold.}
 \begin{description}
\item i) {\rm {If a form $F \in \I'$ \ is s.i. in $\I'$, then also $DF\in \I$ is s.i. in $\I$.}}

\item ii) {\rm {$G\in \I''$ \ is s.i. iff  $G\in \I_t, t>p$ \ and $G$ is s.i. as an element of $\I$.}}
\end{description}

\end{Theorem}

\medskip \noindent
{\bf Proof.} i) \ Let $F\in \I'_u, u\leq p-d$, where $d$ is the degree of $D$, be s.i.. Then $F\notin\I'_{u-1}S_1$, because it is an element of a standard basis of $\I'$. As a consequence $FD \in \I_{d+u}, \  FD \notin \I_{d+u-1}S_1$, so that $FD$ can be an element of some standard basis of $\I$. Let

$\basis =(DF_1,...,DF_m, G_1,...,G_n)$ be any basis of $\I$ such that $ F=F_i$. As $ (F_1,...,F_m)$
is a standard basis of $\I'$, we have
$$(\exists t) \ F\M^t \subset (F_1,...,F_{i-1}, F_{i+1},...,F_m).$$
So, for some $t$,  the relation
$$ (DF)\M^t \subset (DF_1,...,DF_{i-1}, DF_{i+1},...,DF_m, G_1,...,G_n)$$
holds.

\smallskip

ii) \ Let $G$ be a s.i. element of $\I''$.
Thanks to Proposition 3.4 of  \cite{BM2} , stating that no element of degree $\alpha(\I)$ can be s.i.,  $G$ cannot be of the form \ $kD,\ k \in K$, \ so that $t= \deg G \geq p+1$. First we observe that, as an element of $\I$, $G$ can belong to a standard basis.
In fact, as it is a form of a standard basis of $\I''$, we have \ $G\notin (\I''_{t-1})S_1 \supseteq (\I_{t-1})S_1$, \ so that $G\notin \I_{t-1}S_1$.
Now, let $\basis =(DF_1,...,DF_m, G_1,...,G_{i-1},G,G_{i+1},...,G_n)$ be any standard basis of $\I$ containing $G$. Then $(D,G_1,...,G_{i-1},G,G_{i+1},...,G_n)$ is a standard basis of $\I''$. The hypothesis of inessentiality of $G$ as an element of $\I''$ implies  that
$$ (\exists t)\ G \M^t \subset (D,G_1,...,G_{i-1},G_{i+1},...,G_n).$$
In other words, for every form $ P\in \M^t$, we have
$$ GP = DV+\sum_{j\not= i} V_jG_j,$$
so that \  $ (V,V_1,...,V_{i-1},-P,V_{i+1},...,V_n) \in Syz \ \I''$.
\ Condition (c) of Theorem 3.1 (\cite {BDM1}) says that $V\in \I'$, so that  $GP\in (DF_1,...,DF_m, G_1,...,G_{i-1},G_{i+1},...,G_n)$; this means that $G$ is s.i. also as an element of $\I$.

Viceversa, let \ $ G \in \I_t, t>p $  \quad  be a s.i. element in $\I$.
 If $\basis = (DF_1,...,DF_m, G_1,...,G_{i-1},G,$ $ G_{i+1},...,G_n)$ is a basis of $\I$ containing $G$, then  $\basis'' =(D,G_1,...,G_{i-1},G,G_{i+1},...,G_n)$ is a basis of $\I''$.
  Thanks to Proposition 5.1 of \cite{BM2}, it is enough to prove that $G$ is inessential with respect to any basis $\tilde{{\cal B}''}= (D,G_1+A_1G,...,G_{i-1}+A_{i-1}G,G,G_{i+1}+A_{i+1}G,...,G_n+A_nG)$ \ for every (degree allowed) choice of $A_1,...,A_{i-1},A_{i+1},...,A_n.$
  As $\tilde{{\cal B}}=(DF_1,...,DF_m,G_1+A_1G,...,G_{i-1}+A_{i-1}G,G,G_{i+1}+A_{i+1}G,...,G_n+A_nG)$ is still a standard basis of $\I$, $G$ is inessential with respect to it.This means that

$(\exists t\in \N) \ G\M^t \subset (DF_1,...,DF_m,G_1+A_1G,...,G_{i-1}+A_{i-1}G,G_{i+1}+A_{i+1}G,...,G_n+A_nG) \subset (D,G_1+A_1G,...,G_{i-1}+A_{i-1}G,G_{i+1}+A_{i+1}G,...,G_n+A_nG).$

As a consequence, $G$ is inessential also with respect to the basis $\tilde{\basis}''$.

\finedim

\begin{Remark} \label{Rem 2.9} {\rm  Theorem \ref{Th 2.8} can also be proved by working on a suitable H.B. matrix of $\I$, taking into account Proposition 1.2 of \cite{BDM2} and Corollary 4.1 of \cite{BM2}.}
\end{Remark}

\begin{Remark} \label{Rem 2.10} {\rm  It may happen that in $\I$ there exist s.i. elements that do not produce s.i. elements in $\I'$ (see Remark \ref {Rem 3.16})}
\end{Remark}

\begin{Remark} \label{Rem 2.11} {\rm  For every $\I \in F[n]$, the maximality condition  required in Theorem \ref{Th 2.8} is verified at any degree.}
\end{Remark}
\bigskip

 Proposition \ref{Prop 2.2} suggests a situation in which all the elements of every basis of $\I$ are essential because the columns of its H.B. matrix are "short", so that they cannot generate a power of $\M$.

 \begin{Corollary} \label{Cor 2.12} {\rm Let $\I$ be a perfect height 2 ideal of $S= K[X_1,...,X_n]$. Each of the following conditions is enough to guaranty that in any standard basis of $\I$ all the elements are essential:
 \begin{description}
 \item{i) $\nu({\I}) < n+1$}
 \item{ii) $\alpha(\I) < n$}
 \item{ iii) $e(\I) < \displaystyle{\frac{n(n+1)}{2}}$}.
 \end{description}}

 \end{Corollary}

 \medskip \noindent
{\bf Proof.}  i) and ii) are the statement of Corollary 5.1 and Remark in \cite{BM2};  iii) comes from the inequality $\displaystyle{\frac{\alpha(\alpha+1)}{2}} \leq e(\I)$ , where $\alpha= \alpha(\I)$, just observing that $e(\I) < \displaystyle{\frac{n(n+1)}{2}}$ implies $\alpha < n$.

\begin{Remark} {\rm It is easy to find examples of ideals with $e(\I) =\displaystyle{ \frac{n(n+1)}{2}}$ containing inessential elements in some standard basis; see, for instance, Example 3.1 in \cite{BM2}, where $n=3, e=6$.}
\end{Remark}

Proposition 3.4 of \cite{BM2} says that in degree $\alpha(\I)$ no element is s.i.. So, the existence of a basis of essential elements is assured if the generators of degree $>\alpha$ are essential. Such a condition is verified when in the degree matrix \ $\partial M (\I) = (d_{ij}), i=1,...,\nu(\I)-1,  j=1,...,\nu(\I)$ \ the inequality $d_{h,\nu(\I,\alpha)} \leq 0$ is verified for $h=\nu(\I)-n$ (and, as a consequence, for $h<\nu(\I)-n$), because it assures that the columns \ $C_j, j\geq \nu(\I,\alpha),$ \  have at most $n-1$ elements different from zero. This justifies the following statement.

\begin{Proposition} \label{Prop 2.13} {\rm  Let $\I$ be a perfect height $2$ ideal of $S=K[X_1,...,X_n]$, with degree matrix \ $\partial M (\I) = (d_{ij}), i=1,...,\nu(\I)-1,  j=1,...,\nu(\I)$. If \ $d_{\nu(\I)-n,\ \nu(\alpha,\I)} \leq 0$, then $\I$ has a basis of essential elements.}
\end{Proposition}

A consequence of Proposition \ref{Prop 2.13} is the following statement.

\begin{Corollary} \label{Cor 2.14} {\rm  Let $\I$ be a perfect height $2$ ideal of $S=K[X_1,...,X_n]$. If $$e(\I) < \frac{n(n+3)}{2},$$ then $\I$ has a standard basis whose elements are all essential.}
\end{Corollary}

\medskip \noindent
{\bf Proof.} Taking into account the inequality \ $ \displaystyle{\frac{\alpha(\alpha+1)}{2}} \leq e(\I)$, \ we see that the hypothesis implies $\alpha \leq n$. In case \ $\alpha < n$ \ we apply Corollary \ref{Cor 2.12} ii). In case $\alpha = n$ and $\nu = \nu(\I)< \alpha + 1$\ we apply  Corollary \ref{Cor 2.12} i). So, the only case to be considered is \ $ \alpha = n, \nu = n+1$. In this situation the degree matrix \ $\partial M(I) = (d_{ij})$ satisfies the conditions $d_{i,i+1}=1,\ i=1,...,n$. Taking into account the rule of computation of $e(\I)$ starting from \ $\partial M(\I)$ (see \cite{G}), it is easy to verify that the only values of  $d_{ii}$ compatible with  the hypothesis are the following ones:

\begin{description}
\item{a) $d_{ii} = 1, \ i=1,...,n$}
\item{  b)  $d_{ii} = 1,\ i\not= i_0, d_{i_0 i_0} = 2$ for some $i_0 \not= 1.$}
\end{description}

\noindent In case a) the ideal is generated in degree $\alpha$, so that we apply Proposition 3.4 of \cite{BM2}.

\noindent In case b) we have necessarily \ $d_{i_0(i_0+1)}=0$, so that Proposition \ref{Prop 2.13} can be used.
\finedim

\begin{Remark}\label{Rem 2.15} {\rm If the inequality of Corollary \ref{Cor 2.14} is not satisfied, there exist examples of ideals with s.i. elements. For instance, let us consider in $S=K[X_1,...,X_n]$ the ideal $\I$, with H.B. matrix

$$  M(\I) = \begin{pmatrix} X^2_2  & -X_1  &    &     &    &  \\
                                    &  X_2  &-X_1&     &    &  \\
                                    &  X_3  & X_2& -X_1    &    &  \\
                             ....   &  .... &....& ....&... & .... \\
                                    &  X_n  &    &     & X_2 & -X_1 \end{pmatrix},$$
where the unwritten entries are zero forms.

\noindent $\I$ satisfies the condition $ e(\I)= \displaystyle{\frac{n^2+3n}{2}}$ and its second generator is s.i..

We observe that the ideals, with multiplicity \ $ e(\I)= \displaystyle{\frac{n^2+3n}{2}}$ , \ that do not admit a basis of essential elements must necessarily have as a degree matrix the one defined by
$$d_{11}=2; \ d_{ii}=1,\ i\not= 1;\  d_{i(i+1)}=1, \ i=1,...,n, $$
so that they have only one generator in degree $\alpha$.

On the other side, it is possible to produce ideals with a basis of essential elements and with no upper limit on $e(\I)$. For instance, every ideal  $\I$ whose \ $ \partial M(\I)$ \ is defined by
$$d_{i(i+1)}=1, i=1,...,n; \ d_{ii}=1,\ i=1,...,n-1;\  d_{nn}=h\geq 2 $$
satisfies the condition of Proposition \ref{Prop 2.13} and has multiplicity $ e(\I)= \displaystyle{\frac{n(n+1)}{2}+h-1}$, which is arbitrarily large if $h \gg 0$.
}\end{Remark}

\bigskip

We see that if \ $ e(\I)\geq \displaystyle{\frac{n^2+3n}{2}}$ \ the situation is hard to be examined , also if $\I$ is a perfect height $2$ ideal. That is a reason why we restrict our attention to the subfamily ${\cal F}[n]$ (see section 2), starting with $n=2$.

\section{ An e-maximal basis of $\I \in {\cal F}[2]$}

Relation (1.10) in Remark 1 to Theorem 17 in \cite{BM1} gives a good description of every $\I\in {\cal F}[2]$. With some change of notation, we rewrite it as follows:

\begin{eqnarray}
 \lefteqn \I \ =\ \Phi_1...\Phi_r S_{\beta_0} S + ... +\Phi_i...\Phi_r S_{\beta_{i-1}} S +\Phi_{i+1}...\Phi_r S_{\beta_i}S+...+ \Phi_r S_{\beta_{r-1}} S +S_{\beta_r}S   =  \nonumber   \\
         \sum_{i=0}^r \Phi_{i+1}...\Phi_{r +1} S_{\beta_i}S,
          \label{eq 3.1}
\end{eqnarray}
where $\Phi_i$ is a form of degree $\delta_i, \ \Phi_{r+1}=1$ and $S_t$ is the subset of $S=K[X,Y]$ consisting of the forms of degree $t$.

Let us denote $\Delta_0=0, \Delta_i= \delta_1+...+ \delta_i, i=1,...,r$ and $\Delta_r=\delta$ the degree of $\Phi = \Phi_1...\Phi_r$, so that we have
$$\delta=\sum_{i=1}^r \delta_i,  \quad   \beta_i = \beta_{i-1} + \delta_i + t_i ,$$
where $t_i = \alpha_i - \alpha_{i-1} > 0, \ i=1,...,r$  is the difference between two successive different degrees of the generators appearing in a standard basis and $\alpha_0 = \alpha(\I) = \alpha$.

In (\ref{eq 3.1}), \  $r+1$ is the number of distinct elements appearing in any degree vector $\textbf{a}$ of a standard basis of $\I$; moreover, we have
$$\textbf{a}=((\beta_0+\delta)^{[\beta_0+1]}, ...,(\beta_i+\delta-\Delta_i)^{[\delta_i]},...,\beta_r^{[\delta_r]}) = ( \alpha_0^{[\beta_0+1]},\alpha_1^{[\delta_1]},
...,\alpha_i^{[\delta_i]},...,\alpha_r^{[\delta_r]}),$$
where $c^{[n]}$ is the sequence $(c,...,c)$, with $c$ repeated $n$ times.

The degree matrix \ $\partial M (\I) = (d_{ij})$ is completely determined by its elements in position $(i,i+1)$  (which are necessarily $1$, as $M(\I)$ is a $\alpha \times (\alpha +1)$ matrix) and by  $\textbf{a}$, or, equivalently, by  its elements in position $(i,i)$, which are

$d_{ii} =1$ \  if \ $i \not= \beta_0 +1 + \Delta_j, j=0,...,r-1$,

$d_{ii}= t_{j+1} +1$ \ if $i=  \beta_0 +1 +\Delta_j.$

\bigskip

Our aim is to produce an e-maximal basis of $\I$ (see Definition \ref{ Def 2.4} ), that allows to prove the following theorem.

 \begin{Theorem} \label{Th 3.1} {\rm  Let $\I$ be as in (\ref{eq 3.1}) and let

 \begin{equation}
 \Phi = \Phi_1...\Phi_r = H_1^{\mu_1}...H_v^{\mu_v}, \quad \sum_{i=1}^v \mu_i = \delta \geq 1
\label{eq 3.2}
\end{equation}

\noindent be a factorization of $\Phi$ as a product of linear forms pairwise linearly independent.

 The number of s.i. elements appearing in every e-maximal basis of $\I$ is $\delta-v$.
\smallskip

 If $\delta =0$, then $\I = S_{\alpha}S$ does not contain s.i. elements.}

\end{Theorem}

\bigskip
 In order to prove Theorem \ref{Th 3.1} it is convenient (and possible) to produce an e-maximal basis $\basisI$ \ satisfying the following condition.
\begin{equation}
   {\rm  There\ is\ a\ basis\ of\ its\ module\ of\ syzygies \
    linking\ only\ couples\ of\ adjacent\ elements .}
\label{eq.(*)}
\end{equation}

  The Hilbert matrix corresponding to such a basis of $ Syz(\basisI)$ will be called \textit{the canonical matrix} of $\basisI$ or  \textit{a canonical matrix }of $\I$.

  Condition (\ref{eq.(*)})  will be of help in checking that $\basisI$ is an e-maximal basis.

 \medskip \noindent
  Let us consider first two special cases, useful to face the general situation.
\medskip

{\it Case 1}.  \  $\I = S_{\alpha} S$.
\smallskip

Thanks to Proposition 3.4 of \cite{BM2}, we know that an ideal generated in minimal degree cannot have s.i. elements. However, in the sequel we need an e-maximal basis, satisfying condition (\ref{eq.(*)}), constructed according to the following Proposition. The notation $\hat{L}$ will always mean that the element $L$ is omitted.

\begin{Proposition} \label{Prop 3.4} {\rm Let $\I = S_{\alpha} S$. If  $\{L_0,...,L_{\alpha}\}$ is a set of linear forms, pairwise linearly independent. Then \ $\basisI =(F_i), i=0,...,\alpha, \ F_i = L_0...\hat{L_i}...L_{\alpha}$, is a standard basis, consisting entirely of essential elements, whose canonical matrix is
$$M = (m_{ij}) , \ i=1,...,\alpha, \ j=1,...,(\alpha+1),
 {\rm where}:  m_{ii} = L_{i-1},\ m_{i (i+1)} = -L_i, \ m_{ij} = 0 {\rm otherwise}.$$

}
\end{Proposition}

\medskip \noindent
{\bf Proof.} \ It is immediate to verify that the $F_i$'s are a set of $\alpha +1$ linearly independent elements of $S_{\alpha}$, so that they are a basis of it as a $K$-space.  The rows of $M$ are syzygies linking adjacent elements; as they are linearly independent, they are a basis of $Syz \ (\basisI)$ ( see Hilbert-Burch Theorem, \cite{N}), so that $M$ is a matrix of syzygies of $\I$. The entries of every column $C_i$ generate a principal ideal $\I_{C_i}$; so, Proposition \ref{Prop 2.2} \ says that all the elements of $\basisI$ are essential.

\finedim

{\it Case 2}.  \ $\I$ is generated in two different degrees and in the lower one there is just one generator, so that
\begin{equation}
\I=\Phi S + S_b S, \  \deg \Phi = \delta = \alpha(\I),\ b = \beta(\I) = \delta+t, \ t>0.
\label{eq 3.3}\end{equation}

\smallskip

Let us consider the decomposition of $\Phi$ as in (\ref{eq 3.2}), with $r=1$.

We prove first the following lemma.

\begin{Lemma} \label{Lemma 3.5} {\rm Let $\Phi = H_1^{\mu_1}...H_v^{\mu_v}$  be any form of degree $\delta$ in $S = K[X,Y]$. The $K$-space $S_b$, $ \ b=\delta +t, t\geq 0,$  admits a decomposition

\begin{equation}
S_b = \Phi S_t \bigoplus T , \quad  T=\bigoplus^v_{i=1} T_i,
\label{eq 3.4} \end{equation}

where a $K$-basis of $T_i$ is the ordered set $\basis_i = (F_{ij}), j=1,...,\mu_i$,\ described as follows:

\begin{equation}
F_{ij} = A_{ij} C_i,
\label{eq 3.5} \end{equation}
with
\begin{equation}
A_{ij} = H^{\mu_i-j}_i H^{\mu_{i+1}}_{i+1}...H^{\mu_v}_v U^{j-1},\ GCD(U,H_h) =1,\ h=1,...,v,
\label{eq 3.6} \end{equation}
and $C_i$  any form of degree $t+\mu_1+...+\mu_{i-1}+1, (\mu_0=0)$, satisfying the relation

\begin{equation}
GCD(C_i, H_i) = 1.
\label{eq 3.7} \end{equation}
}\end{Lemma}

\medskip \noindent
{\bf Proof.} \ We use induction on $v$.

For $v=1$  we have \ $\Phi = H^{\mu}, \delta = \mu$\ and the statement becomes
$ T= T_1$, with basis $\basis_1 = (F_{1j}),\ j=1,...,\mu$, where

\begin{equation}
F_{1j} = F_j = A_jC = H^{\mu-j} U^{j-1} C, \ \deg C = t+1,\ GCD (C,H) = 1.
\label{eq 3.8} \end{equation}

It is immediate to prove that $F_1,...,F_{\mu}$ are linearly independent, so we only have to show that $\Phi S_t \bigcap T = (0)$.  For $\Lambda \in S_t$, let us suppose $\Lambda \Phi =\sum_{j=1}^{\mu} a_jF_j = (\sum_{j=1}^{\mu} a_j H^{\mu-j} U^{j-1}) C$. This implies that $H^{\mu}$ must divide $ A= \sum_{j=1}^{\mu}a_j H^{\mu-j} U^{j-1}$. For degree reason, $A$ must be zero, so that $a_j=0, \ j=1,...,\mu$.

Let us suppose the statement true until $v-1$ and prove it for $v$.
We set $\Phi = \Psi H_v^{\mu_v}$ and use the decomposition of case  $v=1$ with $H^{\mu}$ replaced by $H^{\mu_v}_v$, so obtaining
$$S_b = H^{\mu_v}_v S_{b-\mu_v} \bigoplus T_v,$$
where \ $ T_v = (F_{v1},...,F_{v \mu_v})$,\ with $F_{vj} = A_{vj}C_v,\ A_{vj} = H_v^{\mu_v-j} U^{j-1}, \ GCD (C_v, H_v) = 1, \deg C_v = b-\mu_v+1$,\ according to (\ref{eq 3.8}).

Using induction, we have $S_{b-\mu_v} = \Psi S_t \bigoplus T',\ T' = \bigoplus  _{i=1}^{v-1} T_i'$,  where $T_i'$ has the basis $(F'_{ij})$ described in the statement of Lemma \ref{Lemma 3.5}, that is
 $F'_{ij} = H^{\mu_i-j}_i H^{\mu_{i+1}}_{i+1}...H^{\mu_{v -1}}_{v-1} U^{j-1}C_i $.
  So, we finally obtain
$$ S_b = H^{\mu_v}_v (\Psi S_t \bigoplus T') \bigoplus T_v  =  \Phi S_t \bigoplus T,$$
where \ $T=\bigoplus H^{\mu_v}_v T' \bigoplus T_v = \bigoplus^{v-1}_{i=1} H^{\mu_v}_v T'_i \bigoplus T_v$. \ It is immediate to check that $(F_{ij})= (H^{\mu_v}_v F'_{ij}),\ j=1,...,\mu_i,$ \ is the required basis of \ $T_i= H^{\mu_v}_v T'_i,\ i=1,...,(v-1).$

\finedim

\begin{Remark} \label{Rem 3.6} {\rm Each space $T_i$ depends on the choice of the form $C_i$, with the link (\ref{eq 3.7}). So, there are infinitely many decompositions of the type described in (\ref{eq 3.4}). Later on, we will use some of them, properly chosen accordingly to the situation.
}\end{Remark}

\begin{Remark} \label{Rem 3.7} {\rm  The basis of $S_b, \ b=\alpha$, used in Proposition \ref{Prop 3.4} is obtained accordingly to Lemma \ref{Lemma 3.5}, with the choice $\Phi = L_0...L_b, \ t=-1, \ C_i = H_1...H_{i-1}$. In this situation, the first summand of (\ref{eq 3.4}) is empty, so that $S_b=T$.}
\end{Remark}

\begin{Proposition} \label{Prop 3.8} {\rm Let us consider the ideal
$$ \I = \Phi S + S_b S, \ b=\delta +t,\ t>0, \quad \Phi = H_1^{\mu_1}...H_v^{\mu_v}, \deg \Phi = \delta.$$

\begin{description}
\item i){ $\I$ has as a standard basis the set
 $$ \basisI = (\Phi, F_{ij}),\ i=1,...,v, \ j=j(i)=1,...,\mu_i,$$
 where:

 \begin{equation}
 F_{ij} = A_{ij}C_i,
 \end{equation} \label{eq 3.9}

 \begin{equation}
  A_{ij} = H_i^{\mu_i-j} H_{i+1}^{\mu_{i+1}}...H_v^{\mu_v} U^{j-1}, \ \ GCD (H_i,U)=1
 \end{equation}\label{eq 3.10}

\begin{equation}
\ C_1=U^{t+1}, \ C_i = H_1H_2...H_{i-1}U^{\nu_i},
 \ \nu_i = t+\mu_1+...+\mu_{i-1}-i+2, \ i>1,
\end{equation}\label{eq 3.11}
}
\item ii){ The basis  $\basisI$ satisfies  condition (\ref{eq.(*)} ). Its canonical matrix $M (\I)$ has as rows the  basis of syzygies $\{ s_{ij}\},\ i=1,...,v, \ j=j(i) =1,...,\mu_i$, with the lexicographic order, where:

    $ s_{11} = (U^{t+1}, -H_1,0,...,0),$

    $ s_{i1} = (0,...,0,H_{i-1}, -H_i,0,...0),\ i=2,...,v,$ \ $-H_i$ in position $\mu_1+\mu_2+...+\mu_{i+2}$,

$s_{ij} = (0,...,0,U,-H_i,0,...,0), \ i=1,...,v,\ j=2,...,\mu_i$ \ , \ $-H_i$ in position $\mu_1+...+\mu_{i-1}+j+1,\ \mu_0 =0$.

     So, $M(\I)$ looks as follows:

$$M(\I)= \begin{pmatrix} U^{t+1}  &   \\
                                   & {\cal A}    \\
                               {\cal O}   &    \end{pmatrix},$$

where  ${\cal A}  =  (a_{ij})$ \ is a square \ $\delta \times \delta$ \  matrix, whose elements different from \ $a_{ii}, a_{(i+1)i}$ are zero, and \ $(a_{11},..., a_{\delta \delta}) = ( [-H_1]^{\mu_1}, [-H_2]^{\mu_2},...,[-H_v]^{\mu_v})$, \  $ a_{(i+1),i} = -a_{ii} $   if \ $a_{ii} \not = a_{jj}, j>i$ \ \ and \ $a_{(i+1),i} = U$\ otherwise.

\item iii) { The essential elements of $\basisI$ are : $ \Phi, F_{(i, \mu_i)}, i= 1,...,v$.  All the other $\delta-v$ elements of $\basisI$ are s.i..}

}\end{description}

}\end{Proposition}

\medskip \noindent
{\bf Proof.} \begin{description}
\item i) {This assertion is an immediate consequence of Lemma \ref{Lemma 3.5}. In fact, thanks to the inequality $\mu_1 +\mu_2+...+\mu_{i-1} \geq  i-1$, we can choose  $C_i = H_1...H_{i-1}U^{\nu_i}$, so obtaining $F_{ij}$ as a basis of the $K$-space $T$ complementary to $\Phi S_t $ in $S_b$.}

\item ii) { The fact that the $\{ s_{ij}\}$'s are syzygies can be verified with an easy direct computation. Moreover, they are clearly linearly independent, of the expected degree and their number $\delta$ is the rank of the module of syzygies, according to Hilbert theorem.  It is easy to verify that the first maximal minor of $M(\I)$ is $\Phi$ and (apart from a sign) the other maximal minors are the $F_{ij}$'s. Using Proposition \ref{Prop 2.2}, we see immediately that the essential columns of $M(\I)$ (that is the columns corresponding to essential elements , see \cite{BM2})  are the first one and the $(\mu_1+\mu_2+...+\mu_i+1)-th, i=1...v$;
so, the essential elements of $\basis$ are $\{\Phi, F_{i,\mu_i}, i=1...v\}$.}

\item iii) {The proof that all the inessential elements are s.i. is a consequence of the following Lemma \ref{Lemma 3.8}, stated in a form sufficiently general to be used later  in a more general situation. In fact, the submatrix $\cal{A}$ appearing in $M(\I)$ satisfies the hypothesis of Lemma \ref{Lemma 3.8} }.

\end{description}

\finedim

It is convenient to generalize the notion of inessential and strongly inessential columns of a matrix, as considered in \cite {BM2}.

\begin{Definition} \label{Def 3.7} {\rm  Let  $A=(a_{ij}$ be a  matrix whose entries are forms of $K[X_1,...,X_n]$ such that $ \deg a_{ij}- \deg a_{i(j+1)}$ is independent from $i$. A column $C_j $ is inessential when the ideal $\I_{C_j}$ generated by its entries contains a power of the irrelevant ideal. $C_j$ is strongly inessential when every column \ $C'_j= \sum_i\lambda_i C_i =(a'_{ij}), \  \lambda_i \in K[X_1,...,X_n],\  \lambda _j=1,\  \deg a'_{ij}= \deg a_{ij}$,\  replacing $C_j$, is still inessential.}
\end{Definition}

\begin{Lemma} \label {Lemma 3.8} {\rm  Let  $A =(a_{ij}), a_{ij} \in K[X,Y],$ be a square $m\times m$ matrix such that $ \deg a_{ij}- \deg a_{i(j+1)}  \leq 0$  is independent from $i$ and satisfying the conditions:

i)  $a_{ij} = 0$  if  $i\not= j, j+1$

ii) $a_{ii}$ is a linear form $L_i$,

iii) $a_{(i+1)i,}$ is any form $G_i$,  such that  $L_j$ is not a factor of $G_i$ if $L_j\not= L_i$ and $G_i$ is a multiple of $L_i$ iff every $L_j, j>i$ is different from $L_i$.

Then the inessential columns of $A$ are s.i.}
\end{Lemma}

\medskip \noindent

{\bf Proof.} \ The inessential column we are considering is of the form $C_j =^t(0,...,0,L_j,G_j,0,...,0)$, where $L_j$ does not divide $G_j$, so that no $L_i$ divides $G_j$. Let us replace such a $C_j$ with  $C'_j = \sum_i C_i, \ \lambda_j = 1$ and prove that $C'_j$ is still inessential.
If $h$ is the first index for which $\lambda_h \not= 0$, we point our attention on the column $C_h$ (clearly, $h\geq j)$. Let us distinguish two possible situations.
\begin{description}
\item i) {$C_h$ \textit{is essential}, so that $C_h = ^t(0,...,0,-L_h, aL_h, 0,...,0)$. In this case we have $C_h \not= C_j$, so that $h<j$ and the entries of $C_h$ must have the same degree of the corresponding entries of $C_j$ (in particular, $a\in K^*$.)

-  If $\lambda_{h+1} \not= 0$, let us consider $C_{h+1} = _t(0,...0,L_{h+1}, G_{h+1},0,...,0)$. The entries of $C'_j$ in position $(h,j)$  and $(h+1,j)$  are respectively \ $c'_{hj}= \lambda_h L_h$ \ and \ $c'_{(h+1) j}= a\lambda_h L_h + \lambda_{h+1} L_{h+1}$,  so that \ $\I_{C'_j} = (L_h, L_{h+1}) = \M$, as $L_h$ is independent from $L_{h+1}$.

-  If \ $\lambda_{h+1} = \lambda_{h+2}=...   \lambda_{h+u-1}=0, \lambda_{h+u}\not=0, u>1$,\ then necessarily \ $h+u \geq j$,\ so that $C'_j$ has as elements \ $c'_{hj} = \lambda_hL_h$ \ and \ $c'_{h+u} = \lambda_{h+u}L_{h+u}, \ \lambda_{h+u}\in K^*$; as a consequence, also in this case \ $\I_{C'_j}=\M$.}

\item ii) {$C_h$ \textit{is inessential}, so that $C_h = ^t(0,...,0,-L_h, aL_h, 0,...,0)$, where $G_h$ is not divisible for $L_q, q=1,...,m.$ ( As a special case, $C_h$ might coincide with $C_j$.) Let us denote $h+u$ the least integer $v$ for which $L_v = L_h$.

    1- If $u=1$, then $c'_{h,j} = \lambda_h L_h, \ c'_{(h+1)j} = \lambda_h G_h + \lambda_{h+1}L_h$, so that \ $\I_{C'_j} \supseteq (L_h, G_h) \supseteq \M^t$, for some  $t\in \N$.

    2- If $u\not= 1$ but $\lambda _{h+1} = 0$, then $c'_{hj} = L_h,\ c'{(h+1)j} = G_h$, so that \ $\I_{C'_j} \supseteq(L_h,G_h)$, as in the previous case.

    3- If $u\not= 1$,\ $\lambda _{h+1} \not= 0$, then \ $c'_{hj} = \lambda_hL_h, \ c'_{(h+1)j} = \lambda_hG_h + \lambda_{h+1} L_{h+1}$,  where  $\lambda_h \in K^*$ (as $h\leq j$).

    If $\lambda_{h+1}$ is such that  $c'_{(h+1)j}$ is not a multiple of $L_h$, we get $C'_j \supseteq \M^t$ , for some $t\in \N$. However, for some choice of $\lambda _{h+1}$ it may happen \  $c'_{(h+1)j} = L_h P$. \ In fact, if \ $G_h = M_1...M_s$\ is a decomposition of $G_h$ into linear factors, there exists $a\in K^*$ such that \ $\lambda_{h+1} = aM_1...M_{s-1}$\  gives \ $c'_{(h+1)j} = M_1...M_{s-1}(\lambda_h M_s + aL_{h+1})$, \ where $\lambda_hM_s + aL_{h+1} =bL_h$,\ as\  $M_s, L_{h+1}$\ are linearly independent linear forms. Let us observe that such a \ $\lambda _{h+1}$ cannot be a multiple of $L_h$, as $G_h$ is not so. If we replace \ $C_h$ with $C^*_h= \lambda C_h + \lambda_{h+1} C_{\lambda+1} C_{h+1}$ \ and consider \ $C'_j = C^*_h +\sum_{i=h+2}^m \lambda_iC_i$, \ we have a situation very similar to the previous one. In fact \ $c*_{hj} = \lambda_hL_h, \ c^*_{(h+1)j} = PL_h, \ c^*_{h+2,j} = \lambda_{h+1} G_{h+1}$, \ so that $G_h$ is replaced with \ $\lambda_{h+1} G_{h+1}$,\ which is not a multiple of $L_h$. Now, we can repeat the same reasoning until when we find either case 2, if $\lambda_i = 0$ for some $i$ with  $h+1<i\leq h+u$,  or case 1, for $i=h+u$.                            }

\end{description}

\finedim

\begin{Remark}\label{Rem 3.8} {\rm The essential generators of \ $\basisI - \{\Phi\}$  are exactly the ones that do not contain as their factors all the linear factors of $\Phi$; more precisely, $F_{i\mu_i}$ does not contain $H_i$, while it contains as factors $H_j, j\not=i$.}
\end{Remark}

In the sequel we will need also bases slightly different from the one produced in Proposition \ref{Prop 3.8}. We introduce them in the following Remarks.

\begin{Remark} \label{Rem 3.9}  {\rm If, in the definition of $F_{ij}, i>k$, $\C_i$  is replaced by  $\tilde{C}_i = H_1...\hat{H}_k...H_{i-1}U^{\nu_i+1}$  (that is, if $H_k$ is replaced with $U$), then $\tilde{\basis}$, obtained from $\basis$ by replacing $F_{ij}$ with $\tilde{F}_{ij} = A_{ij}\tilde{C}_i$, is still a standard basis, whose Hilbert matrix $\tilde{M}(\I)$ differs from the $M(\I)$ described in Proposition \ref{Prop 3.8} just in the column corresponding to $\tilde{F}_{k \mu_k}$, which becomes $\tilde{C}_k = ^t(0,...,-H_k,U,0,...0)$. The consequence is that the generator $\tilde{F}_{k\mu_k} = F_{k\mu_k}$  now is inessential, while the other generators are changed but remain with unchanged  nature.  $\tilde{\basis}$ is not an e-maximal basis, but it will erase in a splitting (see Remark \ref{Rem 3.16}).}
\end{Remark}

\begin{Remark} \label{Rem 3.10} {\rm Let us observe that the $F_{ij}$'s have $U^{t+1}$ as a common factor. If we replace $U^{t+1}$ by any form $\eta$, of degree $t+1$, such that $G.C.D.(\eta,\Phi) = 1$, the matrix  $M^*(\I)$ corresponding to the new basis $\basis ^*$ differs from $M(\I)$ only in the first column. In particular, $\basis^*$ is still an e-maximal basis.}
\end{Remark}
\medskip

\begin{Remark} \label{Rem 3.11} {\rm Let us produce other H.B. canonical matrices of $\I$, relative to standard bases different from the one described in Proposition \ref{Prop 3.8} . They are defined as follows.

$$M'(\I)= \begin{pmatrix} U^{t+1}  &   \\
                                   & {\cal A'}    \\
                               {\cal O}   &    \end{pmatrix},$$

\noindent  where  ${\cal A'}  =  (a'_{ij})$ \ is a square \ $\delta \times \delta$ \  matrix, whose elements different from \ $a'_{ii}, a'_{(i+1)i}$ are zero, and

 -  $(a'_{11},..., a'_{\delta \delta}) =(-H_{\sigma(1)},...,-H_{\sigma(\delta)})$, \ with $\sigma$ any permutation of the sequence $( [1]^{\mu_1}, ( [2]^{\mu_2},..., ( [v]^{\mu_v})$,

 - $ a'_{(i+1),i} = -a'_{ii} $   if \ $a'_{ii} \not = a'_{jj}, j>i$ \ \ and \ $a'_{(i+1),i} = U$\ otherwise.

 In fact, Lemma \ref{Lemma 3.8} guaranties that all the inessential columns of $M'(\I)$ are s.i. and it is a matter of computation to check that the maximal minors of the new matrix are still the basis of a subspace $T$ such that $\Phi S_t \bigoplus T = S_b$. The maximal minors of $M'(\I)$, different from $\Phi$, apart from a sign are: $({\cal B}'_i= U^{t+1}G_1...G_i\hat{H}_{\sigma(i)} H_{\sigma(i+1)}... H_{\sigma(\delta)}),\ i=1,...,\delta.$ A reasoning analogous to the one in the proof of Lemma \ref{Lemma 3.5} shows that they are linearly independent. In fact  the relation

 $\lambda_1 H_{\sigma(1)}...H_{\sigma(\delta)} + \sum_{i=2}^{\delta-1} = 0, \  (\lambda_1,...,\lambda_{\delta}) \not= (0,...,0)$

 \noindent implies that $G_1$ divides $ H_{\sigma(1)}...H_{\sigma(\delta)}$, against the hypothesis.

 Moreover, let us denote $T'$ the $K$-space generated by $({\cal B}'_1,...,{\cal B}'_{\delta})$. Then $\Phi S_t \bigcap T' = (0)$, because $\Lambda \Phi = \sum a_i U^{t+1}G_1...G_i\hat{H}_{\sigma(i)} H_{\sigma(i+1)}... H_{\sigma(\delta)}),\ \Lambda \not = 0$, \ implies that $U$ must divide $\Phi$, for degree reason, against the hypothesis.

}\end{Remark}

\medskip

\begin{Example} \label{Ex 3.11} {\rm Let us consider the ideal
$$ \I = (H_1^3 H_2^2 H_3)S + S_8 S,$$
where $H_1,H_2,H_3$ are linearly independent linear forms. The basis  considered in Proposition \ref{Prop 3.8} is $\basisI \ = (\Phi, F_{11},  F_{12},  F_{13}, F_{21}, F_{22}, F_{31})$, where
$\Phi = H_1^3H_2^2H_3, \ F_{11} = H_1^2H_2^2H_3U^3,\  F_{12} = H_1H_2^2H_3U^4,\  F_{13} = H_2^2H_3U^5,\  F_{21}= H_1H_2H_3U^5, \ F_{22} = H_1H_3U^6, \  F_{31}=H_1H_2U^6$.

The corresponding H.B. matrix is

$$  M(\I) = \begin{pmatrix} U^3   & -H_1  &    &     &     &   &   \\
                                    &  U   &-H_1&     &     &   &   \\
                                    &      & U  &-H_1 &     &   &    \\
                                    &      &    & H_1 &-H_2 &   &    \\
                                    &      &    &     & U   & -H_2&  \\
                                    &      &    &     &     &  H_2 & -H_3  \end{pmatrix},$$
where the unwritten entries are zero forms.

The essential elements are: $ \Phi, F_{13},F_{22},F_{31}$. All the other elements are s.i..

- If in each generator of degree $8$ we replace $U^3$ by any degree $3$ form $\eta$, with $G.C.D.(\eta,\Phi) = 1$, we obtain a new e-maximal basis.

- If we replace $F_{21},F_{22},F_{31}$ respectively by $ \tilde{F}_{21}= H_2H_3U^6, \tilde{F}_{22}= H_3U^7,\tilde{F}_{31}= H_2U^7$, then in the new matrix $\tilde{M}(\I)$  the $H_1$ in $(4,4)$ position is replaced by $U$. As a consequence, $\tilde{F}_{11}=F_{11}$ and $\tilde{F}_{12}=F_{12}$ are s.i., while $\tilde{F}_{13}=F_{13}$ is inessential, but not strongly and $\tilde{F}_{21}\not=F_{21}$ is s.i..

- If we replace $F_{31}$ by $\tilde{F}_{31} = H_1U^7$ (or, equivalently, in $M(\I)$ the form $H_2$ in position $(6,6)$ is replaced by $U$), then $F_{22}$ becomes inessential (but not strongly), while the nature of the other generators does not change.}
\end{Example}

\bigskip
The two special cases just considered suggest us to afford the general case pointing our attention on the H.B. matrix, more than on the standard basis. We need a decomposition of the $\Phi$'s appearing in (\ref{eq 3.1}) into pairwise independent linear forms, as follows

\begin{equation}
\Phi_k = H^{\mu_{k1}}_{k1}...H^{\mu_{k2}}_{k2}...H^{\mu_{kv_k}}_{kv_k},\ k=1,...,r.
\label{eq 3.12}\end{equation}

Moreover, it is convenient to choose a set of linear forms $\{ U, L_0,...,L_{\beta_0}\}$ such that the elements of the set \ $\{U,L_i,H_{kj}\}, \ i=0,...,\beta_0,\ k=1,...,r,\ j=1,...,v_k$ \ are pairwise linearly independent and define

\begin{equation}
\Phi_0 = L_0...L_{\beta_0}.
\label{eq 3.13}\end{equation}

With this notation we can state the following proposition.

\begin{Proposition} \label{Prop 3.12} {\rm A canonical matrix of the ideal $\I$ of (\ref{eq 3.1}) is the following one

  $$  M(\I) = \begin{pmatrix}  {\cal B} & {\cal O}  \\
                                {\cal C} & {\cal A}  \\ \end{pmatrix},$$

where:
\begin{description}
\item i) ${\cal B}  \in S ^{\beta_0\times(\beta_0+1)}, \ {\cal A} \in S^{\delta\times\delta}, \ {\cal C} \in S^{\delta\times(\beta_0+1)}$, \ ${\cal O}$ is a zero matrix, whose elements are of degree $\leq 0$.

\item ii) ${\cal B}= (b_{ij})$, where:
$b_{ii} = L_{i-1}; \ b_{i(i+1)} = -L_i; \  b_{ij} =0$ if $ j\not=i,i+1$.

\item iii) ${\cal C}=(c_{ij}$, where: $c_{1(\beta_0+1)} = U^{t_1+1}; \ c_{ij}=0$ if $(ij)\not=(1(\beta_0+1))$

\item iv)  $ {\cal A} =(a_{ij})$, where:

- \ $a_{ij}=0$ if $j\not=i, j\not=i-1$,

- \ $(a_{11},...,a_{\delta,\delta}) = ([-H_{11}]^{[\mu_{11}]}, ([-H_{12}]^{[\mu_{12}]},...,([-H_{1v_1}]^{[\mu_{1v_1}]},...,[-H_{k1}]^{[\mu_{k1}]}, ([-H_{k2}]^{[\mu_{k2}]},...,$

\hskip 1cm    $ ([-H_{kv_k}]^{[\mu_{kv_k}]},...,[-H_{r1}]^{[\mu_{r1}]}, ([-H_{r2}]^{[\mu_{r2}]},...,([-H_{rv_r}]^{[\mu_{rv_r}]}),$

- \ $a_{(i+1)i} =  \ -a_{ii}$  \ if \ $a_{ii} \not= a_{jj}, j>i, \ i\not= \Delta_k$,

   \quad $a_{(i+1)i} = -a_{ii} U^{t_k}$  \ if\ $a_{ii} \not= a_{jj}, j>i, \ i= \Delta_k$,

   \quad $a_{(i+1)i} = - U$  \ if\ $(\exists j>i) \ a_{ii} = a_{jj}, \ j\not= \Delta_k, \ k<r$,

    \quad $a_{(i+1)i} =  U^{t_k +1}$  \ if\ $ (\exists j>i) \ a_{ii} = a_{jj},\   i=\Delta_k, \ k<r.$

    \end {description}

\noindent Moreover, the inessential columns of $M(\I)$ are s.i..

}\end{Proposition}

\medskip \noindent
{\bf Proof.} We first observe that the degree matrix of $M(\I)$ is the expected one.

Let us denote $\I$ the ideal generated by the maximal minors of $M(\I)$ and prove that $\I$ is the one described in (\ref{eq 3.1}).

As ${\cal B}$ is the matrix considered in Proposition \ref{Prop 3.4}, it is immediate to see that $\I_{\alpha_0} = \Phi S_{\beta_0}$ and that the minors of ${\cal B}$ are linearly independent.

The minors in degree $\alpha_1$ have as a common factor $\Phi/\Phi_1$. So, it is enough to prove that, divided by their common factor, they are a basis of a subspace $T_1$ of $S_{\beta_1}$ such that $S_{\beta_1}= \Phi_1S_{\beta_1-\delta_1} \bigoplus T_1$. But we are in the situation described in Lemma \ref{Lemma 3.5}, where:

  $t=t_1,\ b=\beta_1, \ \Phi = \Phi_1, \ H_i=H_{1i},$

   $C_i= \Phi_0 L_{\beta_0}^{-1}  U^{t_1+1} a_{21}...a_{(j+1)j}, \ j=\mu_{11}+\mu_{12}+...+\mu_{1(i-1)},\ i=1,...,v_1$.

So, let us suppose the statement true until the degree $\alpha_{k-1}$ and prove it for $\alpha_k$.  Just as in the case $k=1$, we see that all the minors have as a common factor $\Phi_{k+1}...\Phi_r = \Phi/\Phi_1...\Phi_k$. So, it is enough to show that, divided by this factor, they are a basis of a subspace
$T_k$ of $S_{\beta_k}$ such that $ S_{\beta_k} = \Phi_k S_{\beta_k-\delta_k} \bigoplus T_k$. We are again in the situation of Lemma \ref{Lemma 3.5}, with:

$t=t_k,\ b=\beta_k,\ \Phi = \Phi_k,\ H_i = H_{ki} $

$C_i = \Phi_0 L_{\beta_0}^{-1} U^{t_1+1} a_{21}...a_{(\Delta_{k-1} +j+1) (\Delta_{k-1} +j)}, \ j=\mu_{k1}+\mu_{k2}+...+\mu_{k(i-1)}, \ i=1,...,v_k$.

Thanks to Proposition \ref{Prop 2.2}, we immediately see that the inessential columns are exactly the ones in which $a_{(i+1)i}$ is not a multiple of $a_{ii}$ or, equivalently, the ones whose element $a_{ii}$ is equal to some $a_{jj}$, with $j>i$. The proof that they are s.i. is a consequence of Lemma \ref{Lemma 3.8}.

\finedim

Extending the notation used in Proposition \ref{Prop 3.8}, we denote the basis linked to the canonical matrix  of Proposition \ref{Prop 3.12} as follows: $\basisI = ({\cal B}^0,{\cal B}^1,...,{\cal B}^k...,{\cal B}^r)$, where ${\cal B}^0= (F_j^0), \ j=0,...,\beta_0, \ {\cal B}^k= (F_{ij}^k), \ k=1,...,r, \  i=1,...,v_k, \ j=1,...,\mu_{ki}$. With this notation we can state the following corollary.

\begin{Corollary} \label{Cor 3.13} {\rm
\begin{description}
\item i) All the elements of ${\cal B}^0$ are essential. The generator $F^k_{i j} \in {\cal B}^k$ is essential iff  $j = \mu_{ki}$ and $H_{ki}$ is not a factor of it.
\item ii) $\basisI$ is an e-maximal basis and the number of its essential elements in degree  bigger then $\alpha(\I)$ is equal to  the number $v$  of the distinct linear factors appearing in a factorization of $\Phi$.
\item iii) In any e-maximal basis the essential generators appearing in degree $\alpha_k$ are as many as the linear factors of $\Phi_k$ that do not divide $ \Phi_{k+1}...\Phi_r$.
\item  iv) $\I$ admits a basis of essential elements iff $\Phi$ is a product of distinct linear factors.    \end{description}

} \end{Corollary}

\medskip\noindent
{\bf Proof}. \begin{description}
\item i) From Proposition \ref{Prop 3.12} iv) we easily see that the essential columns of ${\cal A}$ are the ones whose entry $a_{hh}$ is different from every  $a_{jj}, j>h$. This happens iff $a_{hh} =-H_{ki}$, where $H_{ki}$ does not appear any more in the diagonal of ${\cal A}$, in position $(j,j), j>h$. A necessary condition for such a situation is that the generator corresponding to that column is of the kind $F^k_{i\mu_{ki}}$. In this case we have: $\prod_{j>h}  a_{jj} = R \Phi_{k+1}...\Phi_r$,\ where $H_{ki}$ is not a factor of $R$. So, the condition characterizing the essential $F^k_{i\mu_{ki}}$'s is that $\Phi_{k+1}...\Phi_r$ is not a multiple of $H_{ki}$. From the equality \ $F^k_{i\mu_{ki}} =  \prod_{j>h}  a_{jj} \prod_{j<h}  a_{(j+1)j}$ we see that the previous condition is equivalent to say that $H_{ki}$ does not divide $F^k_{i\mu_{ki}}$.
 \item ii)  $\basisI$ is an e-maximal basis, because its inessential elements are s.i. (Theorem \ref{Th 2.5}). Moreover, the $H_{ki}$ appearing in an  essential column corresponding to $F^k_{i\mu_{ki}}$ is a linear factor of $\Phi$, making there its last appearing as an element of the diagonal of ${\cal A}$. So, the essential columns of ${\cal A}$ are as many as the distinct linear factors of $\Phi$.
\item iii)  As the number of essential elements in an e-maximal basis does not depend on the e-maximal basis chosen, it is enough to verify the statement on the basis $\basisI$ of Proposition \ref{Prop 3.12}. In the proof of i) we observed that the essential elements of $\basis^k$ are as many as the linear factors  $H_{ki}$ of $\Phi_k$ that are not divisors of $ \Phi_{k+1}...\Phi_r$.
\item iv) This is an obvious consequence of ii).

\end{description}
\finedim

Corollary \ref{Cor 3.13}  completes the proof of Theorem \ref{Th 3.1}.

\begin{Corollary} \label{Cor 3.14} {\rm Let $\I$ be represented as in (\ref{eq 3.1}), with $ \Phi_k= H^{\mu_{k1}}_{k1}...H^{\mu_{kv_k}}_{kv_k}$. If $\tau_k$ is the number of distinct linear factors that $\Phi_k$ has in common with $\Phi_{k+1}...\Phi_r$, then any e-maximal basis of $\I$ has exactly  $\sum_{j=1}^{v_k} (\mu_{kj}-1) + \tau_k$ strongly inessential generators in degree $\alpha_k$.}
\end{Corollary}

Corollary \ref{Cor 3.14} implies that it is possible to find $\I \in F[2]$ with a prescribed number of strongly inessential elements in a prescribed number of sufficiently high degree, as we see in the following proposition.

\begin{Proposition} \label{Prop 3.15} {\rm  Let $(d_1<d_2<...<d_s)$ and $(r_1,r_2,...,r_s)$ be sequences of natural numbers. There exist ideals $\I \in F[2]$ with exactly $r_i$ s.i. elements in degree  $d_i,\ i=1,...,s,$ iff
\begin{equation}
d_1>\sum_{i=1}^s r_i+1.  \label{eq 3.14}
\end{equation}
} \end{Proposition}

\medskip \noindent
{\bf Proof.} \ Let us observe that the minimal degree $\delta$ of a form $\Phi$ satisfying the condition $\delta-v = \sum_{i=1}^s r_i$ is obtained with $v=1$, so that $\Phi$ looks as $\Phi = H^{m+1}$, where $m=\sum_{i=1}^s r_i$ and $H$ is any linear form. So, condition (\ref{eq 3.14}) is necessary. It is also sufficient, because the ideal
\begin{equation}
\I = H^{m+1} S + H^{m+1-r_1} S_{d_1-(m+1-r_1)} +...+ H^{m+1-\sum_{i=1}^jr_i} S_{d_j-(m+1-\sum_{i=1}^jr_i)}+...+S_{d_s}
 \label{eq 3.15}  \end{equation}
obtained with the choice $\Phi_j = H^{r_j},  j=1,...,(s-1), \ \Phi_s = H^{r_s+1}$, \ satisfies the required condition. If $d_1=\sum_{i=1}^s r_i +2$, then (\ref{eq 3.15}) is the unique ideal satisfying the condition. If
$d_1 >\sum_{i=1}^s r_i +2$, there are many other possibilities. In fact, the set of the ideals satisfying the required condition increases with the degree $\delta = v+ m $, or, equivalently, with the number $v$ of different linear factors of $\Phi$. Let us observe that the degree vector of the ideal $\I$ considered in (\ref{eq 3.15}) is the least compatible with the required condition.
\finedim

\begin{Remark} \label{Rem 3.16} {\rm Every $\I \in F[2]$ satisfies condition (\ref{eq 2.1}) ( maximality with respect to Dubreil-Campanella inequality) in each degree $\alpha_i$. So, for every $j$, $\I$ splits into two ideals, $\I' = (\I:(\Phi_{j+1}...\Phi_r))$  and  $\I'' = (\I,\Phi_{j+1}...\Phi_r)$, both elements of $F[2]$. The first $\beta_0+1+\Delta_j$ rows and $\beta_0+\Delta_j$ columns of the matrix $M(\I)$  produced in Proposition \ref{Prop 3.12} \ form a H.B. matrix of $\I'$, whose inessential columns are not necessarily s.i.. In fact, it may happen that a linear factor of $\Phi_i, i\leq j$ does not divide $\Phi_{i+1}...\Phi_j$  but divides $\Phi_{j+1}...\Phi_r$; so the assertion of Remark \ref{Rem 2.10} is justified.
} \end{Remark}

\begin{Examples} \label{Ex 3.17} {\rm In the following examples $U,H,K,L_0,L_1,L_2$ are linear forms, pairwise linearly independent.
\begin{description}

\item {\bf 1-}  \quad  $\I = H^3 K^2 S_2 S + K^2 S_6 S + S_{10} S$.

In this case we have: $\Phi = H^3K^2, \Phi_1 = H^3, \Phi_2 = K^2,\ GCD (\Phi_1,\Phi_2) = 1$. According to Proposition \ref{Prop 3.12}, we get

$$  M(\I) = \begin{pmatrix}   L_0  & -L_1  &         &     &   &  &    &   \\
                                    &  L_1  &-L_2     &     &   &   &    &  \\
                                    &      &     L_2 U &-H     &   &    &  & \\
                                    &      &          & U &-H   &   &      &  \\
                                    &      &         &    & U  & -H    &   &  \\
                                    &      &         &     &   & H U^2 & -K &  \\
                                    &      &        &     &   &        & U & -K  \end{pmatrix}.$$

The corresponding canonical basis is

$\basisI = (H^3K^2(L_1L_2, L_0L_2, L_0L_1); K^2L_0L_1L_2U(H^2, HU,U^2); L_0L_1L_2HU^5(K,U))$.

 There are $3$ s.i. generators, according to the  fact that $\delta=5, v=2$. Let us observe that in this example a s.i. generator gives rise to a s.i. generator in any splitting.

 \item {\rm 2-} \quad $\I = H^3K^2S_2S + HKS_6S + S_{10}S$.

 In this case we have: $ \Phi = H^3K^2, \Phi_1 = H^2K, \Phi_2 = HK$, so that all the linear factors of $\Phi_1$ are also divisors of $\Phi_2$. According to Proposition \ref{Prop 3.12}, we get

 $$  M(\I) = \begin{pmatrix}   L_0  & -L_1  &    &     &     &   &  &    &   \\
                                    &  L_1  &-L_2&     &     &   &   &    &  \\
                                    &      &     L_2 U &-H     &   &    &  & \\
                                    &      &          & U &-H   &   &      &  \\
                                    &      &         &    & U  & -K    &   &  \\
                                    &      &         &     &   &  U^3 & -H&  \\
                                    &      &        &     &   &        & H & -K  \end{pmatrix}.$$

The corresponding canonical basis is

$\basisI= (H^3K^2(L_1L_2, L_0L_2, L_0L_1);H K L_0 L_1 L_2 U( H K , U K, U^2); L_0L_1L_2U^6(K,H))$.

  There are $3$ s.i. generators, according to the  fact that $\delta=5, v=2$.

Let us observe that in this case the splitting in degree $p=8$ gives rise to the  matrices

 $$  M(\I') = \begin{pmatrix}   L_0  & -L_1  &    &     &     &   &     \\
                                    &  L_1  &-L_2&     &     &   &     \\
                                    &          &L_2 U &-H     &   &     \\
                                    &           &     & U &-H   &     \\
                                    &            &   &    & U  & -K      \\
                                     \end{pmatrix},$$

 and, respectively,

$$  M(\I'') = \begin{pmatrix}  U^9 & -H  &  \\
                                    & H  &-K  \\
                                      \end{pmatrix}.$$

 The $6$-th column of $M(\I')$ is inessential, but not s.i.. A H.B. matrix of $\I'$, whose corresponding inessential generators are s.i., can be obtained from $M(\I')$ just by replacing $U$ with $H$ in its $6$-th column.

 \item {\rm 3-} \quad $\I = H^3S_2 + H^2S_4 +H S_8 + S_{10}$.

 In this case we have: $\Phi = H^3, \Phi_1= \Phi_2 = \Phi_3 = H$. According to Proposition \ref{Prop 3.12}, we get:
 $$  M(\I) = \begin{pmatrix}   L_0  & -L_1  &    &     &     &    &    \\
                                    &  L_1  &-L_2&     &     &      &  \\
                                    &          &L_2 U &-H     &      & \\
                                    &           &     & U^4 &-H       & \\
                                    &         &    &   &     U^2 & -H
                                     \end{pmatrix}.$$

There are $2$ s.i. generators and the corresponding canonical basis is:

$\basisI= ( H^3( L_1L_2, L_0L_2, L_0L_1); H^2L_0L_1L_2U; HL_0L_1L_2U^5; L_0L_1L_2U^7)$

The splittings in degrees $8$ and $9$ produce a new inessential, but not s.i., element in $\I'$.

\end{description}

}\end{Examples}

\section{ Behaviour of $\I \in F[3]$  with respect to essentiality: special cases and examples}

According to Theorem 1.5 of \cite{BM1}, every element $\I\in F[3]$ has a shape very similar to the one described in (\ref{eq 3.1}) for the elements of $F[2]$. The difference is that $S_{\beta_i}$ is replaced by a linear subspace $T_{\beta_i}\subseteq S_{\beta_i}$ of  $S= K[X,Y,Z]$, where  dim $T_{\beta_i} = \beta_i+1$. The subspaces $T_{\beta_i}$ are characterized by Theorem 3.4 of \cite{BM1}. That theorem says that, up to a change of cohordinates, every element $\I\in F[3]$ is generated by the maximal minors of an $\alpha \times (\alpha+1)$ matrix, obtained by lifting to $K[X,Y,Z]$ a H.B. matrix of its image $\bar{\I} \subset K[X,Y]=\bar{S}$, modulo a regular linear form Z. So, $\I$ looks like
\begin{eqnarray}
 \I \ =\ \Phi_1...\Phi_r T_{\beta_0} S + ... +\Phi_i...\Phi_r T_{\beta_{i-1}} S +\Phi_{i+1}...\Phi_r T_{\beta_i}S+...+ \Phi_r T_{\beta_{r-1}} S +T_{\beta_r}S   =  \nonumber   \\
         \sum_{i=0}^r \Phi_{i+1}...\Phi_{r +1} T_{\beta_i}S, ,
          \label{eq 4.1}
\end{eqnarray}
where $\Phi_i$  is a form in $S$ and  $\Phi_{r +1}=1$,
and its image modulo $Z$ becomes
\begin{eqnarray}
  \bar{\I}  =\  \bar{\Phi}_1...\bar{\Phi}_r \bar{S}_{\beta_0}\bar{S} + ... +\bar{\Phi}_i...\bar{\Phi}_r \bar{S}_{\beta_{i-1}} \bar{S} +\bar{\Phi}_{i+1}...\bar{\Phi}_r \bar{S}_{\beta_i}\bar{S}+...+ \bar{\Phi}_r \bar{S}_{\beta_{r-1}} \bar{S} +\bar{S}_{\beta_r}\bar{S}   =  \nonumber   \\
         \sum_{i=0}^r \bar{\Phi}_{i+1}...\bar{\Phi}_{r +1}
         \bar{S}_{\beta_i} \bar{S},\ \bar{\Phi}_{r +1}=1 .
          \label{eq 4.2}
\end{eqnarray}

The problem of stating if a lifting of \ $\bar{\I} \in F[2]$  to $\I \in F[3]$  preserves the strong inessentiality of the entries of $\basis (\bar {\I})$ becomes a lifting problem of H.B. matrices, which seems not easy to be solved. So, we start to consider a very special case. More precisely, we focus our attention on the ideals $\I\in F[3]$  with the largest number of s.i. generators in any e-maximal basis. If $\alpha= \alpha(\I) =\a(\bar{\I})$ is the minimal degree of the generators of $\I$, we will see that the maximal expected number is  $\alpha-2$; we'll prove that such a number is reached. Let us first state a property for every homogeneous saturated ideal $\I$ of $S=K[X_1,...,X_n]$.

\begin{Proposition} \label{Prop 4.1} {\rm Let $\basis =(b_1,...,b_h, c_1,...,c_k), \ k\geq 1$, be an e-maximal basis of the saturated homogeneous ideal $\I\subset S = K[X_1,...,X_n] $, where $b_1,...,b_h$ are essential and $c_1,...,c_k$  are s.i. elements. The condition depth $\I = r$  implies $h>r$.}
\end{Proposition}

\medskip\noindent
{\bf Proof}. Thanks to Corollary 5.3 of \cite{BM2},  $(c_1,...,c_k)$ is an inessential set (\cite{BM2}, Def. 5.2), so that $\I=(b_1,...,b_h)^{sat}$. As  depth $\I$ = depth $(b_1,...,b_h)$, the hypothesis implies $h\geq r$; however, the equality holds iff $(b_1,...,b_h)$ is a c.i. and, as a consequence, a saturated ideal, against the hypothesis $k\geq 1$.

\finedim

Choosing $h=2$, we get immediately the following statement.

\begin{Corollary} \label{Cor 4.2} {\rm The largest possible number of s.i. generators in an e-maximal basis of an ideal $\I\in F[3]$ is \ $\alpha(\I) - 2$.}
\end{Corollary}

If $\a(\I) = 2$, then $\I \in F[3]$ has $3$ generators and Corollary \ref{Cor 4.2} says that in every e-maximal basis they are essential. Let us point our attention on the case $\a(\I)>2$.

\medskip
\noindent  We will use the following  \underbar{Notation}
$$ {\cal S} =  \{ \I \in F[3] : \nu_e(\I) = 3, \ \a(\I) >2 \},$$  where $\nu_e(\I)$ denotes the number of essential elements of any e-maximal basis of $\I$ (see \cite{BM2} Def. 5.1).

\medskip
Let us observe that any dehomogenization $\I_*$ with respect to a regular linear form of an ideal $\I\in {\cal S}$ has just $3$ generators (that is the least number for a non complete intersection), while the number of generators of $\I$ is the maximum allowed by Dubreil's inequality.

\medskip

\begin{Proposition} \label{ Prop 4.3} {\rm  For every  $\I \in {\cal S}$ the form $\Phi$ appearing in  (\ref {eq 4.1}) is of degree $\delta = \alpha (\I)$  and $\Phi$  has necessarily one of the following shapes:
\begin{description}
\item {i) $\Phi = H^{\delta}$,}
\item { ii)  $\Phi = H^rK^s, \ r+s = \delta,$}
\item {iii) $\Phi = C^{\gamma}, 2\gamma = \delta$,}
    \end{description}
where $H$ and $K$ are independent linear forms and $C$ is a quadratic irreducible form in $K[X,Y,Z]$.}
\end{Proposition}

\medskip
\noindent {\bf Proof.}  \ An immediate consequence of Proposition \ref{Prop 2.7} is that if $\I\subset K[X,Y,Z]$ has $\alpha -2$  s.i. generators, then the number of s.i. generators of its quotient $\bar{\I}$ modulo a regular linear form is either $\alpha -2$ or $\alpha -1$. Applying Theorem \ref{Th 3.1} to $\bar{\I}$, we immediately get
\begin{equation}
\delta-v = \alpha -1, \ \delta \leq \alpha,    \label{eq 4.3}
\end{equation} or
\begin{equation}
  \delta-v = \alpha -2, \  \delta \leq \alpha . \label{eq 4.4}
\end{equation}

Relation (\ref{eq 4.3}) is equivalent to $\delta=\alpha, \ v=1$, while relation (\ref{eq  4.4}) gives two possible situations:

\begin{equation}
\delta = \alpha , \ v=2 \label{eq 4.5}
\end{equation}
and
\begin{equation}
  \delta = \a-1,\ v=1. \label{eq 4.6}
\end{equation}

Let us verify that (\ref{eq 4.6}) is not realized. In fact in this case we have \  $(\bar{\I})_{\alpha} = \bar{H}^{\alpha-1}\bar{S}_1$  and the s.i. generators lie all in degree bigger than $\a$; so, a splitting in degree $\a$ ( see Theorem \ref{Th 2.8}) gives rise to an ideal $\I''$, with $\a (\I'') = \a - 1$ and $\a-2$ s.i. generators, against Corollary \ref{Cor 4.2}.

So, $\Phi$ must be a form, of degree $\alpha$, whose quotient modulo any regular linear form splits into a product of powers of at most two different linear factors. This means that the curve $\Phi = 0$  meets a generic  line in at most two different points, so that $\Phi$ is necessarily as described in i), ii), iii).
\finedim

\begin{Remark}\label{Rem 4.1}

\qquad 1- {\rm We do not have examples in which the situation iii) appears. Let us observe that it requires every $\nu(\I,j), \ j=1,...,r$, to be a power of $2$.}

2- {\rm Proposition \ref{ Prop 4.3} says that the schemes corresponding to ideals with $\alpha -2$ s.i. generators lie necessarily either on a multiple line or on two multiple lines or (may be) on a multiple irreducible conic. However this condition is not sufficient. For instance, by lifting the canonical matrix $M$ of an  ideal $\J$ of $K[X,Y]$  with $\alpha-1$ s.i. generators with $M$ itself, we obtain a basis for an ideal $\I \subset K[X,Y,Z]$ without inessential elements.}

3- {\rm In case i), $\bar{\Phi}$ has the same structure of $\Phi$, for every regular form $L$, while in case ii) and iii)  we have generically $\bar{\Phi} = \bar{H}^r \bar{K}^s$, but we get $\bar{\Phi} = \bar{M}^{r+s}$, $M$ a linear form in $K[X,Y]$, iff the line $L=0$ either is tangent to $\Phi=0$ (case iii))  or meets it in its singular point (case ii)). So, it is possible to represent an element $\I \in {\cal S}$ as a lifting of ideals $\bar{\I}\subset k[X,Y]$ such that $ \bar{\Phi} = \bar{M}^{\delta}$, except for the case in which  $\Phi=H^rK^s$ and the intersection between the two lines $H=0$ and $K=0$ is in the support of the corresponding scheme.}
\end{Remark}

\begin{Proposition} \label{Prop 4.4}{\rm If $\I$ is an element of ${\cal S}$, let us consider its splitting into $\I'$ and $\I''$ (see Theorem \ref{Th 1.2}), in degree $\alpha_k $, with  $k<r$   if  $\deg \Phi_r \geq 2$ and  $k< r-1$  if $\deg \ \Phi_r = 1$. Then $\I''$ is still an element of ${\cal S}$.}
\end{Proposition}\
\medskip \noindent
{\bf Proof.}  Theorem \ref{Th 2.8} says that if $\basisI = (\Phi, \basis^1, \basis^2,...,\basis^r)$ \ is an e-maximal basis of $\I$, then $\basis ({\I''}) =(\Phi'',\basis^{k+1},...,\basis^r)$, where $\Phi'' = \Phi_{k+1}...\Phi_r$, is an e-maximal basis of $\I''$ and the forms of $\basis(\I'')$, different from $\Phi''$, maintain the same nature they had in $\basisI$.
So, as $\basisI$  has two essential elements different from $\Phi$, $\basis(\I'')$ cannot contain more then two essential elements different from $\Phi''$; the hypothesis on the choice of $k$ guaranties that it has at least $3$ elements; then, Corollary \ref{Cor 4.2} says that $\basis (\I'')$ must have exactly $3$ essential elements, so that $\I'' \in {\cal S}$.
\finedim

\begin{Corollary} \label{Cor 4.5} {\rm In any e-maximal basis of an ideal $\I\in {\cal S}$\ the degree of the three essential generators are: $\a =\a(\I),\ \a_r(\I)$  and either $\a_r(\I)$  or $\a_{r-1}(\I)$. The latter possibility takes place iff in degree $\a_r(\I)$ there is just one generator.}
\end{Corollary}
\medskip \noindent
{\bf Proof.}  It is enough to apply Proposition \ref{Prop 4.4}, with $k=r-1$ if in degree $\a_r$\ the basis $\basisI$ \ contains at least two forms and with \ $k=r-2$\ otherwise. In fact the two essential elements of $\basis(\I'')$,  not in minimal degree, must be essential also in $\basisI$.
\finedim

Now, let us produce examples of ideals of ${\cal S}$.

\begin{Proposition} \label{Prop 4.5} {\rm Every $\bar{\I} \subset K[X,Y]$  with $\alpha-1$  s.i. generators in any e-maximal basis has at least a lifting in ${\cal S}$.}
\end{Proposition}

\medskip \noindent
{\bf Proof.} \ After a possible change of coordinates, a H.B. matrix of $\bar{\I}$ is

\begin{equation}
  M(\bar{\I}) = \begin{pmatrix}   Y^{t_0}  & -X  &    &     &     &   &  &    &   \\
                                    &      Y^{t_1} &-X&     &     &   &   &    &  \\
                                    &      &     Y^{t_2} &-X     &   &    &  & \\

                                    &...   &...      & ...   & ...  & ...    &...   &  \\
                                    &      &         &     &   &  Y^{t_{\delta-2}} & -X&  \\
                                    &      &        &     &   &        & Y^{t_{\delta-1}} & -X  \end{pmatrix}                 \label{ eq 4.8} \end {equation}

A lifting of $M(\bar{\I})$ with $\alpha-2$ \  s.i. columns is

\begin{equation}
  M(\I) = \begin{pmatrix}   Y^{t_0}  & -X  &    &     &     &   &  &    &   \\
                                  ZP_1  &      Y^{t_1} &-X&     &     &   &   &    &  \\
                                  ZP_2  & Z^{t_1+t_2-1}     &     Y^{t_2} &-X     &   &    &  & \\
                                  ZP_3  &   &Z^{t_2+t_3-1}      & Y^{t_3}   & -X  &     &   &  \\
                                               &...   &...      & ...   & ...  & ...    &...   &  \\
                                  ZP_{\delta-2}  &      &         &     & Z^{t_{\delta-3}+t_{\delta -2}-1 } &  Y^{t_{\delta-2}} & -X&  \\
                                  ZP_{\delta-1}  &      &        &     &   & Z^{t_{\delta-2}+t_{\delta-1}-1}       & Y^{t_{\delta-1}} & -X  \end{pmatrix}                 \label{ eq 4.8} \end {equation}

                      \smallskip            \noindent
where $P_i\in K[X,Y,Z]$ is a form of degree $\sum _{j=0}^i t_j -i-1$ and, as usual, the unwritten entries are zero. $M(\I)$ is obtained from $M(\bar{\I})$ by leaving unchanged the last two columns and replacing the zero entries in position $(i+1,i), i=2,...,\delta-1,$ with $Z^{t_{i-1}+ t_i-1}$ and the ones in position $(i,1)$ with $ZP_{i-1}$. It is immediate to verify that the second column is s.i. and Lemma 5.3 of \cite{BM2} says that the same reasoning can be repeated for the following ones with three non zero entries.

\finedim

Finally, let us point our attention on the ideals $\I\in {\cal S}$ with the smallest $\alpha(\I)$  allowing the presence of s.i. generators. Corollary \ref{Cor 4.2} implies that if $\I$ has a s.i. generator then $\alpha(\I) \geq 3$. So, let us look for the ideals with $\alpha(\I) =3$ and just one s.i. generator in every e-maximal basis; they are the  elements of ${\cal S}$ with the smallest number of generators. Let us consider the special case of generators in two different degrees. Proposition \ref{ Prop 4.3} says that, apart from a coordinates change, they can be obtained by lifting an ideal of one of the following types
\begin{equation}
\bar{\I}_1 = X^3S +S_{\beta}S, \ \beta >3,\end{equation} \label{eq 4.12}

\begin{equation}
\bar{\I}_2 = X^2YS +S_{\beta}S,\ \beta >3.\end{equation} \label{eq 4.13}

Let us first consider all the required liftings of $\bar{\I}_1$ or, equivalently, all the liftings  $M(\I_1)$ of the matrix

\begin{equation}
  M(\bar{\I_1}) = \begin{pmatrix}   Y ^t  & -X  &  0  & 0     \\
                                  0  &    Y &-X&  0         \\
                                  0  &  0    &     Y &  -X   \\
                                                         \end{pmatrix}, \quad t=\b-2, \label{eq 4.14} \end{equation}
having a s.i. column. $M(\I)$ has the following shape ( see \cite{BM1})

\begin{equation}
  M(\I_1) = \begin{pmatrix}   Y ^t +ZP_1(X,Y,Z) & -X + a_{11}Z &  a_{12}Z & a_{13}Z     \\
                                  ZP_2(X,Y,Z)  &    Y + a_{21}Z &-X + a_{22}Z &  a_{23}Z         \\
                                  ZP_3(X,Y,Z) &  a_{31}Z   &  Y+ a_{32}Z &  -X +a_{33}Z \\
      \end{pmatrix}, \label{eq 4.15} \end{equation}
where $a_{ij}\in K,\ \deg P_i = t-1$. The forms $P_1,P_2,P_3$ can be chosen arbitrarily among the ones of degree $t-1$, so that we just have to characterize the matrix ${\cal A} = (a_{ij}), \ i,j=1,2,3$. As the first column of $M(\I_1)$ is essential for every choice of the $a_{ij}$'s, let us consider the second and third columns.
 The second column is s.i. iff the forms
 $$ -X + a_{11}Z + \l_2 a_{12}Z +\l_3 a_{13}Z,\quad Y + a_{21}Z + \l_2(-X+a_{22}Z) + \l_3 a_{23}Z, \quad  a_{31}Z + \l_2(Y+a_{23}Z) + \l_3(-X+a_{33}Z)$$
 are linearly independent, for every choice of $ \l_2, \l_3$ or, equivalently, iff the matrix
$$ {\cal B} = \begin{pmatrix} -1  &  0  & a_{11}+\l_2a_{12}+\l_3a_{13}  \\
                           -\l_2  &  1 & a_{21}+\l_2a_{22}+\l_3a_{23}  \\
                           -\l_3  &  \l_2 & a_{31}+\l_2a_{32}+\l_3a_{33} \\  \end{pmatrix}$$
has determinant different from zero. Such a condition gives the relation
\begin{eqnarray}
 -a_{12}\l_2^3 - a_{13}\l_2^2 \l_3 + (a_{22}-a_{11})\l_2^2 + a_{13}\l_3^2 + (a_{12} +a_{23})\l_2 \l_3 + \nonumber \\
  (-a_{32} + a_{21})\l_2 + ( -a_{33} + a_{11}) \l_3 - a_{31} \not= 0.\end{eqnarray} \label{eq 4.16}
 An easy computation shows that the matrices  ${\cal A}$ for which this condition is satisfied are
 \begin{equation}
{\cal A} = \begin{pmatrix} a_{11}  &  0  &  0  \\
                            a_{21}  &  a_{11} & 0 \\
                               a_{31}  &  a_{21}  & a_{11} \end{pmatrix}, \qquad a_{31} \not= 0 . \label{eq 4.16} \end{equation}
Considerations very similar to the previous ones lead to the conclusion that the second column is s.i. iff the matrix ${\cal A}$ has the following shape

 \begin{equation}
{\cal A} = \begin{pmatrix} a_{11}  &  a_{12}  &  0  \\
                            a_{21}  &  a_{11} & -a_{12} \\
                              0     &  a_{21}  & a_{11} \end{pmatrix}, \qquad a_{12} \not= 0 . \label{eq 4.17}
\end{equation}

In case (\ref{eq 4.16}), after the coordinate change  $-a_{11}Z +X =X', \ a_{21} Z+Y = Y',\ a_{31}Z = Z' $ and dropping the apostrophes, the required matrix can be written
\begin{equation}
  M(\I_{11}) = \begin{pmatrix}   Y ^t+ Z Q_1 & -X  &  0  & 0     \\
                                  Z Q_2  &    Y &-X &  0         \\
                                  Z Q_3  &  Z    &     Y &  -X   \\
                                                         \end{pmatrix} \label{eq 4.18} \end{equation}

In case (\ref{eq 4.17}), after the coordinate change  $-a_{11}Z +X =X', \ -a_{21} Z+Y = Y',\ a_{12}Z = Z' $ and dropping the apostrophes, the required matrix can be written
\begin{equation}
  M(\I_{12}) = \begin{pmatrix}   Y ^t+ Z Q_1 & -X  &  Z  & 0     \\
                                  Z Q_2  &    Y &-X &  -Z        \\
                                  Z Q_3  &  0    &     Y &  -X   \\
                                                         \end{pmatrix}.\quad t=\b-2, \label{eq 4.19} \end{equation}

 Let us observe that both schemes relative to $\I_{11}$ \ and $\I_{12}$  are supported at at most \ $t+1$\ points lying on a triple line ($X = 0$ in our basis) and that their multiplicity is \ $ e(\I) = 3+3t$ (\cite{G}).

 \smallskip
 With a very similar computation it is possible to see that, apart from a coordinate change, a lifting $\I_2$ of $\bar{\I}_2$ belongs to ${\cal S}$ iff its H.B. matrix has the following shape
\begin{equation}
  M(\I_{2}) = \begin{pmatrix}   (X+Y) ^t+ Z Q_1 & -X  &  0  & 0     \\
                                  Z Q_2  &    X+Y &-X &  0         \\
                                  Z Q_3  &  Z    &     X &  -Y   \\
                                                         \end{pmatrix}, \label{eq 4.20} \end{equation}
where $P_1,P_2,P_3$ are forms of degree $t-1$ in $K[X,Y,Z]$.
The corresponding schemes, still of multiplicity $e(\I) =3t+3$, are all supported at two different lines ($X=0$ and $Y=0$). The intersection of the two lines is one of the points in the support of the scheme; as a consequence, the ideals cannot be obtained by lifting an ideal of type (\ref{eq 4.14}).

\medskip
The characterization of the elements of ${\cal S}$ with $\a(\I)>3$ is more difficult to be faced, also for ideals generated in two degrees. In fact, the request of (\ref{eq 4.14})  (and the analogous for the lifting of $\bar{\I}_2$) are replaced by the requirement that a system of non linear equations $\{ E_u(a_{ij}, \l_v)=0 \}$, in a set $\{\l_v, v=1,...,\a -1\}$ of variables,  admits no solutions. Such a condition defines the entries  $a_{ij}$'s of the matrix ${\cal A}$ as the  elements for which the ideal generated by the $E_u$'s is the whole ring $K[\l_1,...,\l_{\a-1}]$.

\end{document}